\documentclass[11pt]{amsart}
\usepackage[T1]{fontenc}
\usepackage{lmodern,amsmath,amssymb,amsthm,mathtools,microtype,comment}
\usepackage[margin=1in]{geometry}
\usepackage[hidelinks]{hyperref}

\numberwithin{equation}{section}
\newtheorem{theorem}{Theorem}[section]
\newtheorem{lemma}[theorem]{Lemma}
\newtheorem{proposition}[theorem]{Proposition}
\newtheorem{corollary}[theorem]{Corollary}
\theoremstyle{remark}
\newtheorem{remark}[theorem]{Remark}

\newcommand{\E}{\mathbb E}
\newcommand{\PP}{\mathbb P}
\newcommand{\R}{\mathbb R}
\newcommand{\CC}{\mathbb C}
\newcommand{\HH}{\mathbb H}
\newcommand{\FF}{\mathcal F}
\newcommand{\Lcal}{\mathcal L}
\newcommand{\Q}{\mathbb Q}
\newcommand{\dist}{\operatorname{dist}}
\newcommand{\diam}{\operatorname{diam}}
\newcommand{\crad}{\operatorname{crad}}
\newcommand{\Cont}{\operatorname{Cont}}
\newcommand{\ind}{\mathbf 1}
\newcommand{\area}{\operatorname{Area}}

\usepackage{mathrsfs} 

\newcommand{\C}{\mathbb C}
\newcommand{\Hh}{\mathbb H}
\newcommand{\N}{\mathbb N}
\newcommand{\Z}{\mathbb Z}

\newcommand{\F}{{\mathcal F}}

\def\eps{\varepsilon}

\def\sem{\setminus}
\def\lin{\overline}

\newcommand{\BGE}{\begin{equation}}
\newcommand{\BGEN}{\begin{equation*}}
\newcommand{\EDE}{\end{equation}}
\newcommand{\EDEN}{\end{equation*}}

\title[Exponential integrability of SLE Minkowski content]
{Exponential Integrability of SLE Minkowski Content}
\author{Dapeng Zhan}
\address{Department of Mathematics, Michigan State University,
East Lansing, MI 48824, USA}
\date{November 12, 2026}

\begin{document}
\begin{abstract}
Let $0<\kappa<8$, $d=1+\kappa/8$, and $q=2/(2-d)$.
Let $\mu$ be the $d$-dimensional Minkowski content measure of either
chordal $\mathrm{SLE}_\kappa$ in $\Hh$ from $0$ to $\infty$ or
two-sided whole-plane $\mathrm{SLE}_\kappa$ through $0$.
We prove that there exist constants $c_\kappa,C_\kappa\in (0,\infty)$ such that
\[
 \E\exp\left\{c_\kappa
 \left(\frac{\mu(D)}{(\diam D)^d}\right)^q\right\}\le C_\kappa
\]
for every deterministic bounded Borel set $D$ of positive diameter in the corresponding domain.
The exponent $q$ in this moment bound is optimal in both settings.
We also show that, almost surely,
$\mu(B(z,r))=O\bigl(r^d(\log(1/r))^{1/q}\bigr)$ as $r\downarrow0$,
locally uniformly in $z$.
\end{abstract}
\maketitle

\section{Introduction}
Fix $0<\kappa<8$ and write
\begin{equation}\label{eq:dq-intro}
 d=1+\frac{\kappa}{8}\in(1,2),\qquad
 p=2-d=1-\frac{\kappa}{8},\qquad
 q=\frac{2}{p}=\frac{2}{2-d}>2.
\end{equation}
The natural parametrization of $\mathrm{SLE}_\kappa$ was introduced by
Lawler and Sheffield via a Doob--Meyer construction~\cite{LS}.
Lawler and Zhou subsequently established it for the full range
$0<\kappa<8$~\cite{LZ}.  Lawler and Rezaei then proved that the
$d$-dimensional Minkowski content of chordal $\mathrm{SLE}_\kappa$
exists and gives this parametrization~\cite{LR}.  Let $\gamma$ be
chordal $\mathrm{SLE}_\kappa$ in $\HH$ from $0$ to $\infty$, and let
$\mu_\gamma$ denote its occupation measure under this parametrization.
The natural length of $\gamma$ in any bounded set has finite moments
of all positive integer orders~\cite{RZ}.  We strengthen this to a
position-uniform stretched-exponential estimate and then prove the
same estimate for the self-similar stationary-increment two-sided
whole-plane process.

\begin{theorem}[Chordal SLE]\label{thm:main}
There are constants $c_\kappa,C_\kappa\in(0,\infty)$ depending only on $\kappa$ such that, for every
deterministic bounded Borel set $D\subset\HH$ with $L=\diam D>0$,
\begin{equation}\label{eq:main}
 \E\exp\left\{c_\kappa
 \left(\frac{\mu_\gamma(D)}{L^d}\right)^q\right\}\le C_\kappa.
\end{equation}
\end{theorem}

Let $\eta:\R\to\CC$ be a two-sided whole-plane $\mathrm{SLE}_\kappa$ from $\infty$ to $\infty$ through $0$, parametrized by Minkowski content and normalized by $\eta(0)=0$.  By \cite[Corollary~4.7]{ZhanLoop}, $\eta$ is self-similar with stationary increments (sssi), with self-similarity index $1/d$.  Together with the
rotational invariance of two-sided whole-plane $\mathrm{SLE}_\kappa$, this
gives, for $a\in\R$ and $z_0\in\C\sem\{0\}$,
\begin{equation}\label{eq:sssi-intro}
 (\eta(a+t)-\eta(a))_{t\in\R}\stackrel{\rm law}= (\eta(t))_{t\in\R},
 \qquad
 (\eta(|z_0|^d t)/z_0)_{t\in\R}\stackrel{\rm law}= (\eta(t))_{t\in\R}.
\end{equation}
We call this process sssi $\mathrm{SLE}_\kappa$.

Let $\mu_\eta$ denote the occupation measure of $\eta$, that is, the
pushforward of Lebesgue measure on $\R$ under $\eta$.
The following theorem is the two-sided whole-plane analogue of
Theorem~\ref{thm:main}.

\begin{theorem}[sssi SLE]\label{thm:sssi}
There are constants $c_\kappa,C_\kappa\in(0,\infty)$ depending only on $\kappa$ such that, for every
deterministic bounded Borel set $D\subset\CC$ with $L=\diam D>0$,
\begin{equation}\label{eq:sssi-main}
 \E\exp\left\{c_\kappa
 \left(\frac{\mu_\eta(D)}{L^d}\right)^q\right\}\le C_\kappa.
\end{equation}
\end{theorem}

For chordal SLE, it was conjectured in \cite[p.~184, remark~4 after Theorem~1.2]{RZ} that for every bounded domain $D\subset\HH$, one has $\E[e^{\lambda\mu_\gamma(D)}]<\infty$ for some $\lambda>0$. Since $q>1$, Theorems~\ref{thm:main} and~\ref{thm:sssi} give finite exponential moments for every $\lambda>0$ for both processes. Thus, for either occupation measure $\mu\in\{\mu_\gamma,\mu_\eta\}$, the moment generating function of $\mu(D)$ extends to an entire
function, and the moments of $\mu(D)$ determine its law.

Related exponential tail estimates for discrete approximations of Minkowski content were recently obtained by Cai, Li, Shiraishi, and Xia~\cite[Theorem~3.4 and Section~7]{CaiLiShiraishiXia}. For ordinary LERW in dimensions two and three, they prove exponential upper tails for the number of boxes of a given scale visited by the path, uniformly in both the lattice mesh and the spatial resolution.  At the critical scaling, these box-counting quantities are discrete analogues of Minkowski content.  Their tail estimate implies, by integration, that the corresponding critically rescaled box count admits a nontrivial exponential moment.  In dimension two, this estimate, together with the natural-parametrization convergence of LERW to radial $\mathrm{SLE}_2$ proved by Lawler and Viklund~\cite{LawlerViklund}, should imply a nontrivial exponential moment for the Minkowski content of radial $\mathrm{SLE}_2$.

The next theorem shows that the exponent $q$ in
Theorems~\ref{thm:main} and~\ref{thm:sssi} is optimal.

\begin{theorem}[Sharpness of the exponent $q$]\label{thm:sharpness}
There exist constants $a,b,C>0$ and an integer $m_0\ge1$, depending only on $\kappa$, such that the following estimates hold for every integer $m\ge m_0$.
\begin{subequations}\label{eq:sharp}
\textup{(i)} For chordal $\mathrm{SLE}_\kappa$, if $z\in\HH$, $y=\Im z$, and $0<r<y$, then
\begin{equation}\label{eq:sharp-chordal}
 \PP\left\{r^{-d}\mu_\gamma(B(z,r))\ge a m^p\right\}
 \ge C\left(\frac{r}{y}\right)^p
       \left(\frac{y}{|z|}\right)^{8/\kappa-1}e^{-bm^2}.
\end{equation}
\textup{(ii)} For sssi $\mathrm{SLE}_\kappa$, if $z\in\CC$ and $r>0$, then
\begin{equation}\label{eq:sharp-sssi}
 \PP\left\{r^{-d}\mu_\eta(B(z,r))\ge a m^p\right\}
 \ge C\left(1\wedge \frac{r}{|z|}\right)^p e^{-bm^2},
\end{equation}
\end{subequations}
where we interpret $r/|z|$ as $+\infty$ when $z=0$.
Consequently, for either occupation measure $\mu$ and every nonempty deterministic open set $D$ in the corresponding domain, there is a constant $c_D\in(0,\infty)$ depending only on $\kappa$ and $D$, such that
\begin{equation}\label{eq:critical-coefficient-finite}
 \E\exp \{c_D {\mu(D)} ^q \}=\infty.
\end{equation}\end{theorem}

For $\mu\in\{\mu_\gamma,\mu_\eta\}$ and every deterministic bounded Borel
set $D$ of positive diameter, the exponential moment bounds in
Theorems~\ref{thm:main} and~\ref{thm:sssi}, together with
\[
x^s\le \Big(\frac{s}{e c_\kappa q}\Big)^{s/q} e^{c_\kappa x^q},
\qquad x\ge0,\quad s>0,
\]
applied to $x=\mu(D)/\diam(D)^d$, imply
\[
\|\mu(D)\|_{L^s}\le C_\kappa (\diam D)^d s^{1/q},
\qquad s\ge1.
\]
This growth rate is sharp.  In the chordal case, fix for instance
$K=B(2i,1)$, while in the sssi case let $K=B(0,1)$.  Taking
$m=\lceil\sqrt{s}\rceil$ in \eqref{eq:sharp} gives
\[
\|\mu(K)\|_{L^s}
\ge a m^p \PP\{\mu(K)\ge a m^p\}^{1/s}
\ge C a m^p e^{-b m^2/s}
\gtrsim s^{p/2}=s^{1/q}.
\]
Consequently, $\|\mu(K)\|_{L^s}\asymp_\kappa s^{1/q}$ as  $ s\to\infty$.

\subsection*{Outline of the proofs}

Section~\ref{sec:chordal} proves Theorem~\ref{thm:main}.  The starting
point is the chordal SLE$_\kappa$ Green martingale $M_t(z)=G_{H_t}(z;\gamma(t),\infty)$.
Freezing its value when $z$ first comes within live Euclidean distance
$\varepsilon$ of the trace gives a spatial integral with a
stretched-exponential tail by the estimate of \cite{OAI}.  The main
bridge is to identify its limit with the Minkowski-content occupation
$\mu_\gamma(K)$ of a fixed interior disk $K$.  We first obtain $L^1$
convergence from convergence of the Doob--Meyer compensators, and then
upgrade this to $L^2$ convergence using the $L^2$ estimates for
approximating Minkowski contents from \cite{LR}.  A position-uniform
Whitney cover then extends the resulting disk estimate to arbitrary
bounded Borel subsets of $\mathbb H$, including sets arbitrarily close
to the real line.

Section~\ref{sec:sssi} deduces Theorem~\ref{thm:sssi} from the chordal
result, rather than developing a separate stopped-density theory for
whole-plane $\mathrm{SLE}_\kappa(2)$.  A conformal localization argument
first transfers Theorem~\ref{thm:main} to interior balls in simply
connected slit domains.  The decomposition of the rooted SLE loop
measure with a point marked according to Minkowski content
\cite[Theorem~4.1(iv)]{ZhanLoop} then yields the required estimate for a
centered ball.  Finally, stationary increments give an abstract
mass-transport identity which extends this estimate to arbitrary
bounded Borel sets.  This last step uses no Markov property.

Section~\ref{sec:sharpness} proves Theorem~\ref{thm:sharpness}.  The
exponent $q=2/(2-d)$ is inherited from the stopped-density estimate in
\cite{OAI}, while the matching lower tail is obtained from the
dense-visit construction of \cite[Proposition~6.4]{OAI}.  Letting the
auxiliary sausage mesh tend to zero transfers the finite-mesh lower
bound to Minkowski-content occupation.
The sssi lower bound then follows by conditioning on the first arm
and applying the chordal lower bound to the second arm in the
remaining domain, together with scaling and rotational invariance.

Finally, Section~\ref{sec:uniform-local} derives uniform local
occupation upper bounds from the exponential moment estimates by a
covering argument and the Borel--Cantelli lemma. Combining these
with known continuity estimates yields two-sided local occupation
bounds.

\section{Exponential integrability for chordal SLE} \label{sec:chordal}
\subsection{Stopped Green densities and the mesh tail}\label{sec:chordal-mesh}
Let $\gamma$ be chordal $\mathrm{SLE}_\kappa$ from $0$ to $\infty$
in $\HH$, in half-plane-capacity time.
Let $H_t$ be the surviving domain, i.e., the unbounded component of $\HH\sem \gamma([0,t])$, and set
\begin{equation}\label{eq:green}
 M_t(z)=G_{H_t}(z;\gamma(t),\infty),
 \qquad
 G_U(z;a,b):=a_\kappa\,\crad_U(z)^{-p}S_U(z;a,b)^{8/\kappa-1},
\end{equation}
with zero continuation on swallowing. Let $f:U\to\HH$ be a conformal map taking the boundary prime ends $a,b$ to $0,\infty$, respectively. Here $\crad_U(z)=2\Im f(z)/|f'(z)|$ denotes the full conformal radius, and we set $S_U(z;a,b):=\sin(\arg f(z))$. Both expressions are independent of the choice of $f$, since any two such maps differ by a positive multiplicative constant. The constant $a_\kappa>0$ is denoted by $\hat c$ in \cite[Section~2.3, (8)]{LR}.

For a planar set $S$ and $r>0$, write $S^r=\{z\in\mathbb C:\dist(z,S)\le r\}$ for its closed Euclidean $r$-neighborhood. For a bounded set $S\subset\CC$, its $d$-dimensional upper Minkowski content is
\[
 \lin{\Cont}_d(S):=\limsup_{r\downarrow0}r^{d-2}\area(S^r).
\]
If the limit exists, we call it the $d$-dimensional Minkowski content of $S$ and denote it by $\Cont_d(S)$.

Let $S\subset\CC$ be locally closed, that is, the intersection of an open set and a closed set. Following \cite[Definition~6.3]{ZhanBoundaryMC}, a Borel measure $M$ on $\CC$ with $M(\CC\setminus S)=0$ is a $d$-dimensional Minkowski content measure of $S$ if, for every compact $K\subset S$, we have $M(K)<\infty$ and, whenever $M(\partial_S K)=0$, $\Cont_d(K)$ exists and equals $M(K)$. Here $\partial_S K$ is the relative boundary of $K$ in $S$. Such a measure, when it exists, is unique \cite[Lemma~6.5]{ZhanBoundaryMC}.

A curve $\xi:I\to\CC$, where $I$ is a real interval, is parametrized by $d$-dimensional Minkowski content if $\Cont_d(\xi([s,t]))=t-s$ for all $s,t\in I$ with $s\le t$. Chordal $\mathrm{SLE}_\kappa$ in $\HH$ from $0$ to $\infty$ almost surely admits this parametrization \cite[Theorem~1.1]{LR}, and its occupation measure $\mu_\gamma$ under this parametrization is the Minkowski content measure of its trace \cite[Example~6.9]{ZhanBoundaryMC}. For the sssi $\mathrm{SLE}_\kappa$ curve $\eta$ of Theorem~\ref{thm:sssi}, the occupation measure $\mu_\eta$ is likewise the $d$-dimensional Minkowski content measure of its trace \cite[Remark~2.4 and Remark~2.4]{ZhanLoop}.  In
Section~\ref{sec:chordal}, we write $\mu=\mu_\gamma$ when there is no risk of confusion.

For $z\in\HH$ and $0<\varepsilon<\Im z$, let
$ \tau_\varepsilon(z):=\inf\{t\ge0:z\in(\gamma[0,t])^\varepsilon\}$
with the convention $\inf\varnothing=\infty$. Note that when $\tau_\eps(z)<\infty$, $|\gamma(\tau_\eps(z))-z|=\eps$ and $z\in H_{\tau_\eps(z)}$. Define
\begin{equation}\label{eq:weights}
 W_\varepsilon(z)=M_{\tau_\varepsilon(z)}(z)
          \ind_{\{\tau_\varepsilon(z)<\infty\}},
 \qquad
 V_\varepsilon(K)=\int_K W_\varepsilon(z)\,dA(z).
\end{equation}
Koebe's theorem gives
\begin{equation}\label{eq:pointwise-cap}
 0\le M_{t\wedge\tau_\varepsilon(z)}(z)\le a_\kappa\varepsilon^{-p},
 \qquad
 0\le W_\varepsilon(z)\le a_\kappa\varepsilon^{-p}.
\end{equation}
By \eqref{eq:pointwise-cap}, the stopped local martingale $M_{t\wedge \tau_\eps(z)}(z)$ is bounded and hence a true martingale. On $\{\tau_\varepsilon(z)=\infty\}$, $M_t(z)\to0$ almost surely as $t\to\infty$ (see, e.g., \cite[Proof of Lemma 3.3]{OAI}), so the bounded martingale has terminal value $W_\varepsilon(z)$, which implies that
\BGE \E[W_\varepsilon(z)|\F_t]=M_{t\wedge \tau_\eps(z)}(z),\quad t\ge 0.\label{WM-condition}\EDE

We call a disk $B(z_0,r)$ admissible if $8r<\Im z_0$. Such a disk is compactly contained in $\HH$.

\begin{proposition}\label{input:tail}
There are $c_0,C_0\in(0,\infty)$, depending only on $\kappa$, such that for
any admissible disk $K=B(z_0,r)$, and any $\eps\in (0,\Im z_0- r)$,
\begin{equation}\label{eq:input-tail}
 \PP\{V_\varepsilon(K)>u r^d\}
 \le C_0\exp(-c_0u^q),\qquad u\ge 0.
\end{equation}
\end{proposition}

Since $\dist(K,\R)=\Im z_0-r>7r$, scaling $K$ to unit radius gives the fixed boundary clearance $\rho_{\rm col}=7$ required in \cite[Section~4, especially (4.3)]{OAI}. Scaling back and multiplying the density by $a_\kappa$, we obtain \eqref{eq:input-tail} from \cite[(4.3)]{OAI} for $0<\varepsilon<7r$ and $u\ge u_0$. For $7r\le\varepsilon<\Im z_0-r$, \eqref{eq:pointwise-cap} gives $V_\varepsilon(K)\le\pi r^2a_\kappa\varepsilon^{-p}\le\pi a_\kappa 7^{-p}r^d$. Thus, taking $u_0\ge\pi a_\kappa 7^{-p}$ makes the event in \eqref{eq:input-tail} empty for $u\ge u_0$ in this range of $\varepsilon$. Finally, increasing $C_0$ so that $C_0e^{-c_0u_0^q}\ge1$ gives the bound for $0\le u<u_0$ as well.

\subsection{Identification with Minkowski content} \label{sec:chordal-identification}

Let $(\FF_t)_{t\ge0}$ be the usual augmentation of the filtration generated by the Loewner driving function. Fix a disk $K\Subset\HH$. Let $\Theta_t(K)$ denote the natural time that $\gamma$ spends in $K$ up to capacity time $t$, and set $X:=\Theta_\infty(K)=\mu(K)$. Since $K$ can be covered by finitely many dyadic squares compactly contained in $\HH$, \cite[Theorem~3.1]{LR} gives $X\in L^2$. The conditional expectation identity in \cite[Theorem~1.1]{LR} is stated in terms of Minkowski content. Its occupation-time form \cite[(4.3)]{ZhanDecomposition} gives
\begin{equation}\label{eq:natural-dm}
 U_t:=\int_K M_t(z)\,dA(z)=\E[X\mid\FF_t]-\Theta_t(K).
\end{equation}
The process $\Theta_t(K)$ is continuous and increasing by \cite[Theorem~1]{LZ}. The spatial integral defining $U_t$ in \eqref{eq:natural-dm} is almost surely continuous in $t$ \cite[Section~1.3, immediately after (7)]{LZ}.

For $0<\varepsilon<\dist(K,\R)$, define the accumulated stopped-density mass
\begin{equation}\label{eq:Aeps}
 A_t^\varepsilon
 =\int_K W_\varepsilon(z)\ind_{\{\tau_\varepsilon(z)\le t\}}\,dA(z),
 \qquad A_\infty^\varepsilon=V_\varepsilon(K),
\end{equation}
and the residual potential
\begin{equation}\label{eq:Ueps}
 U_t^\varepsilon
 =\int_K M_t(z)\ind_{\{t<\tau_\varepsilon(z)\}}\,dA(z).
\end{equation}
The integrands in \eqref{eq:Aeps} and \eqref{eq:Ueps} are jointly measurable in $(t,z,\omega)$, so $A_t^\varepsilon$ and $U_t^\varepsilon$ are jointly measurable in $(t,\omega)$.

By \eqref{eq:weights} and \eqref{eq:pointwise-cap}, $V_\varepsilon(K)$ is bounded by $a_\kappa\varepsilon^{-p}\area(K)$ and is therefore in $L^2$. For each $t\ge0$, \eqref{WM-condition} and conditional Fubini give
\BGE\label{eq:eps-dm}
 N_t^\varepsilon:=\int_K M_{t\wedge\tau_\varepsilon(z)}(z)\,dA(z)
 =A_t^\varepsilon+U_t^\varepsilon
 =\E[V_\varepsilon(K)\mid\FF_t].
\EDE
Taking expectations at $t=0$ and using \eqref{eq:natural-dm} gives
\[
 \E V_\varepsilon(K)=\int_K G_{\HH}(z;0,\infty)\,dA(z)=\E X.
\]

The process $A^\varepsilon$ is adapted and increasing. It is continuous by bounded convergence, since $W_\varepsilon$ is bounded on $K$ and, for each finite $t$, the set $\{z\in K:\tau_\varepsilon(z)=t\}$ lies on the circle $\{z:|z-\gamma(t)|=\varepsilon\}$, which has planar area zero. Since $\varepsilon<\dist(K,\R)$, we have $A_0^\varepsilon=0$. The stopped integrand in \eqref{eq:eps-dm} is continuous in $t$ and uniformly bounded by \eqref{eq:pointwise-cap}, so bounded convergence gives continuity of $N^\varepsilon$. Hence $U^\varepsilon$ is continuous as well.

We use the continuous case of the compensator convergence theorem of Dellacherie and Meyer \cite[Chapter~VII, \S20]{DM}, as stated in \cite[Theorem~D]{LZ}.

\begin{proposition} \label{input:compensators}
Let $U^n$, $n\ge1$, and $U$ be continuous nonnegative supermartingales on a filtered probability space satisfying the usual conditions. Suppose that $U$ is of class $D$ and that, almost surely, $U_t^n\uparrow U_t$ for every $t\ge0$. Let $A^n$ and $A$ denote their respective Doob--Meyer compensators, normalized to start at zero. Then $A_T^n\to A_T$ in $L^1$ for every almost surely finite stopping time $T$.
\end{proposition}

Directly from \eqref{eq:natural-dm} and \eqref{eq:Ueps},
\[
 U_t-U_t^\varepsilon
 =\int_K M_t(z)\ind_{\{\tau_\varepsilon(z)\le t\}}\,dA(z)\ge0.
\]
For every finite $t$, as $\varepsilon\downarrow0$, the integrand in \eqref{eq:Ueps} increases to $M_t(z)$ at each $z\in K$. Indeed, if $z\in H_t$, then $\tau_\varepsilon(z)>t$ whenever $\varepsilon<\dist(z,\gamma[0,t])$, while if $z\notin H_t$, then $M_t(z)=0$ by convention. Hence, by monotone convergence,
\begin{equation}\label{eq:monotone}
 0\le U_t^\varepsilon\uparrow U_t
 \quad\text{as }\varepsilon\downarrow0,\quad\text{for every }t\ge0.
\end{equation}

\begin{theorem}[Stopped-density identification]\label{thm:bridge}
For every fixed disk $K\Subset\HH$,
\begin{equation}\label{eq:bridge}
 V_\varepsilon(K)\longrightarrow\mu(K)\quad\text{in }L^2
 \quad\text{as }\varepsilon\downarrow0.
\end{equation}
\end{theorem}

\begin{proof}
Fix any sequence $\varepsilon_n\downarrow0$ with $\varepsilon_1<\dist(K,\R)$. By \eqref{eq:eps-dm} and \eqref{eq:natural-dm}, $A^{\varepsilon_n}$ and $\Theta(K)$ are the compensators of $U^{\varepsilon_n}$ and $U$, respectively. Up to the normalization of the Green function, $U_t$ is the process $\Psi_t(K)$ defined after \cite[(7)]{LZ}. Thus \cite[(16)]{LZ} gives $\sup_T\E[U_T^2]<\infty$, where $T$ ranges over all almost surely finite stopping times, so $U$ is of class $D$. Thus \eqref{eq:monotone} and Proposition~\ref{input:compensators} give $A_T^{\varepsilon_n}\to\Theta_T(K)$ in $L^1$ for every almost surely finite stopping time $T$.

Taking expectations in \eqref{eq:eps-dm} and \eqref{eq:natural-dm} and using $U_t^{\varepsilon_n}\le U_t$, we obtain
\[
 \E[V_{\varepsilon_n}(K)-A_t^{\varepsilon_n}]
 =\E U_t^{\varepsilon_n}
 \le\E U_t
 =\E[X-\Theta_t(K)].
\]
Since $A_t^{\varepsilon_n}\le V_{\varepsilon_n}(K)$ and $\Theta_t(K)\le X$, the triangle inequality and the preceding estimate give
\begin{align*}
 \E\bigl[\lvert V_{\varepsilon_n}(K)-X\rvert\bigr]
 &\le \E\bigl[\lvert A_t^{\varepsilon_n}-\Theta_t(K)\rvert\bigr]
      +\E[V_{\varepsilon_n}(K)-A_t^{\varepsilon_n}]+\E[X-\Theta_t(K)]\\
 &\le \E\bigl[\lvert A_t^{\varepsilon_n}-\Theta_t(K)\rvert\bigr]
      +2\E[X-\Theta_t(K)].
\end{align*}
Since $0\le\Theta_t(K)\uparrow X$ and $X\in L^1$, the second term on the right tends to zero as $t\to\infty$. For each fixed $t$, taking $T\equiv t$ in Proposition~\ref{input:compensators} gives $A_t^{\varepsilon_n}\to\Theta_t(K)$ in $L^1$, so the first term on the right tends to zero as $n\to\infty$. So we obtain $V_{\varepsilon_n}(K)\to X$ in $L^1$.

To obtain $L^2$ convergence, cover $K$ by finitely many dyadic squares $Q_1,\ldots,Q_m$ compactly contained in $\HH$. Since $W_\varepsilon(z)=0$ unless $z\in\gamma^\eps$, \eqref{eq:pointwise-cap} gives
\[
 0\le V_\varepsilon(K)
 \le a_\kappa\sum_{j=1}^m \varepsilon^{-p}
 \area(Q_j\cap \gamma^\eps).
\]
By \cite[Theorem~3.1]{LR}, each normalized area on the right converges in $L^2$ as $\varepsilon\downarrow0$. Hence the finite sum converges in $L^2$, so its squares along $(\varepsilon_n)$ are uniformly integrable. The displayed bound gives uniform integrability of $V_{\varepsilon_n}(K)^2$. Since $V_{\varepsilon_n}(K)\to X$ in probability, Vitali's theorem gives $V_{\varepsilon_n}(K)\to X=\mu(K)$ in $L^2$. The sequence was arbitrary, so \eqref{eq:bridge} follows.
\end{proof}

\subsection{From disks to bounded sets} \label{sec:chordal-bounded}

\begin{proposition}\label{prop:local}
Let $c_0$ and $C_0$ be as in Proposition~\ref{input:tail}. Then for every admissible disk $K=B(z_0,r)$,
\begin{equation}\label{eq:local}
 \E\exp\left\{\frac{c_0}2\left(\frac{\mu(K)}{r^d}\right)^q\right\}\le 1+C_0.
\end{equation}
\end{proposition}

\begin{proof}
Let $a=c_0/2$. Integration of the tail and \eqref{eq:input-tail} give
\begin{align*}
 \E\exp\{a[V_\varepsilon(K)/r^d]^q\}
 &=1+a\int_0^\infty q u^{q-1}e^{au^q}
       \PP\{V_\varepsilon(K)/r^d>u\}\,du\\
 &\le 1+aC_0\int_0^\infty q u^{q-1}e^{(a-c_0)u^q}\,du
 =1+C_0.
\end{align*}
By the $L^2$ convergence in Theorem~\ref{thm:bridge}, there exists a sequence $\varepsilon_n\downarrow0$ such that $V_{\varepsilon_n}(K)\to\mu(K)$ almost surely. Fatou's lemma along this sequence gives \eqref{eq:local}.
\end{proof}

\begin{lemma}\label{lem:cover}
Let $D\subset\HH$ have finite positive diameter $L$. Then $D$ can be covered by admissible disks $K_j=B(z_j,r_j)$ with
\begin{equation}\label{eq:cover-budget}
 \sum_j r_j^d\le C_d L^d,\qquad C_d:=\frac{8^{2-d}}{1-2^{1-d}}.
\end{equation}
\end{lemma}

\begin{proof}
We call a square with side length $l$ and center $z_0$ admissible if $8l< \Im z_0$. Such a square is contained in the admissible disk $B(z_0,l)$. Thus, it suffices to show that $D$ can be covered by admissible squares with side length $r_j$ such that \eqref{eq:cover-budget} holds.

 The horizontal and vertical spans of $D$ are at most $L$. Choose a closed  real interval $I$ of length $L$ containing the  projection of $D$ to $\R$, and let $h=\inf_{z\in D}\Im z$. In each of the grids below, we cover every square by the concentric open disk of radius equal to its side length.

If $h\ge L$, then $D\subset I\times[h,h+L]$. The square $I\times[h,h+L]$ can be covered by $64$ small squares of side length $r_j=L/8$. Each small square is admissible since its center $z_j$ satisfies $\Im z_j >h\ge L=8r_j$. In addition, we have \eqref{eq:cover-budget} since $\sum_j r_j^d=64 (L/8)^d\le C_dL^d$.

If $h<L$, then $D\subset I\times(0,2L]=\bigcup_{k=0}^\infty I\times [2^{-k} L,2^{1-k}L]$. For $k\ge0$,  divide $I\times[2^{-k} L,2^{1-k}L]$ into $8$ rows and $2^{k+3}$ columns, giving $2^{k+6}$ squares of side length $2^{-k-3}L$. Each square is admissible since the height of its center is greater than $2^{-k}L=8 r_j$. Summing across layers gives
\[
 \sum_j r_j^d
 =\sum_{k=0}^\infty 2^{k+6} (2^{-k-3}L)^d
 =8^{2-d}L^d\sum_{k=0}^\infty 2^{-k(d-1)}
 =C_dL^d.
\]
\end{proof}

\begin{proof}[Proof of Theorem~\ref{thm:main}]
By Lemma~\ref{lem:cover}, choose the cover of $D$ by admissible disks $K_j=B(z_j,r_j)$ satisfying \eqref{eq:cover-budget}, and put $S=\sum_j r_j^d$, $w_j=\frac{r_j^d}{S}$, and $Y_j=\frac{\mu(K_j)}{r_j^d}$.
Then $\sum_jw_j=1$, $\mu(D)\le S\sum_jw_jY_j$, and
Proposition~\ref{prop:local} gives $\E e^{\frac{c_0}2Y_j^q}\le 1+ C_0$.
Since $q>1$, the function $f(x)=e^{\frac{c_0}2 x^q}$ is increasing and convex on $[0,\infty)$. By Jensen's inequality and Tonelli's theorem,
\[
 \E\exp\left\{\frac{c_0}2\left(\frac{\mu(D)}{S}\right)^q\right\}
 \le \E f\left(\sum_j w_jY_j\right)
 \le \sum_j w_j\E f(Y_j)
 \le 1+C_0.
\]
Since $S\le C_d(\diam D)^d$ by \eqref{eq:cover-budget}, this proves \eqref{eq:main} with $c_\kappa=\frac{c_0}2 C_d^{-q}$ and $C_\kappa=1+C_0$.
\end{proof}

\begin{remark}
Let $D\subset\HH$ be a deterministic bounded Borel set of positive diameter. Since $\mu$ is the $d$-dimensional Minkowski content measure of $\gamma$, \cite[Lemma~6.6(i)]{ZhanBoundaryMC} and the transience of $\gamma$ give
\BGE \lin{\Cont}_d(K\cap \gamma)\le \mu(K\cap\gamma)=\mu(K),\quad\text{for all compact }K\subset\C.\label{lin-Cont}\EDE
This together with $\mu(\R)=0$ implies  $ \lin{\Cont}_d(\gamma\cap D)  \le \mu(\overline D\cap \HH)$.
Since  $\diam(\overline D\cap\HH)=\diam D$, applying \eqref{eq:main} to $\overline D\cap\HH$ gives the same exponential moment bound for $\lin{\Cont}_d(\gamma\cap D)$.
\end{remark}

\section{Exponential integrability for sssi SLE} \label{sec:sssi}

\subsection{Conformal localization and the SLE loop measure}\label{sec:sssi-localization}

We use two properties of Minkowski content measures to transfer Theorem~\ref{thm:main} to the interior of a simply connected domain. If $M$ is the $d$-dimensional Minkowski content measure of $S$ and $G\subset\CC$ is open, then $M|_G$ is the Minkowski content measure of $S\cap G$ \cite[Theorem~6.12]{ZhanBoundaryMC}. If $S\cap G$ is relatively closed in $G$ and $f:G\to\CC$ is conformal and injective, then $f_*(|f'|^d\,M|_G)$ is the Minkowski content measure of $f(S\cap G)$ \cite[Theorem~6.13]{ZhanBoundaryMC}.

Let $\gamma_U$ be chordal $\mathrm{SLE}_\kappa$ in a simply connected domain $U\subsetneq\CC$ from a boundary prime end $a$ to a distinct boundary prime end $b$. Choose a conformal map $f:\HH\to U$ taking $0$ to $a$ and $\infty$ to $b$, and couple $\gamma_U$ with chordal $\mathrm{SLE}_\kappa$ in $\HH$ from $0$ to $\infty$, denoted by $\gamma$, so that $\gamma_U\cap U=f(\gamma\cap\HH)$. By the restriction property \cite[Theorem~6.12]{ZhanBoundaryMC}, $\mu_\gamma|_{\HH}$ is the $d$-dimensional Minkowski content measure of $\gamma\cap\HH$. Since $\gamma\cap\HH$ is relatively closed in $\HH$ by transience, conformal covariance \cite[Theorem~6.13]{ZhanBoundaryMC} gives the $d$-dimensional Minkowski content measure of $\gamma_U\cap U$ as
\begin{equation}\label{eq:interior-mc}
 \mu_U:=f_*\bigl(|f'|^d\,\mu_\gamma|_{\HH}\bigr).
\end{equation}

\begin{lemma}\label{lem:interior-transport}
There are constants $c,C\in(0,\infty)$ depending only on $\kappa$ such that, for every $U,a,b$ as above and every disk $B(z_0,r)$ with $B(z_0,4r)\subset U$,
\begin{equation}\label{eq:interior-transport}
 \E\exp\left\{c\left(\frac{\mu_U(B(z_0,r))}{r^d}\right)^q\right\}\le C.
\end{equation}
\end{lemma}

\begin{proof}
Let $f:\HH\to U$ be as above, put $g=f^{-1}$, and set $B=B(z_0,r)$ and $E=g(B)$. Applying the Koebe distortion theorem to the univalent map $\zeta\mapsto g(z_0+4r\zeta)$ on the unit disk gives
\[
 \frac{\sup_B|g'|}{\inf_B|g'|}
 \le \left(\frac{1+1/4}{1-1/4}\right)^4
 =\left(\frac53\right)^4<8.
\]
Since $\diam E\le 2r\sup_B|g'|$ and $\sup_E|f'|=(\inf_B|g'|)^{-1}$, we obtain
\begin{equation}\label{eq:distortion}
 (\sup_E|f'|)\diam E<16r.
\end{equation}
By \eqref{eq:interior-mc} and \eqref{eq:distortion},
\[
 \frac{\mu_U(B)}{r^d}
 =\frac{1}{r^d}\int_E |f'(w)|^d\,\mu_\gamma(dw)
 \le 16^d\,\frac{\mu_\gamma(E)}{(\diam E)^d}.
\]
Applying Theorem~\ref{thm:main} to $E$ gives \eqref{eq:interior-transport} with $c=16^{-dq}c_\kappa$ and $C=C_\kappa$.
\end{proof}

Let $\Lcal=\mu_\infty^1$ be the $\mathrm{SLE}_\kappa$ loop measure rooted at infinity, with the normalization of \cite[Theorem~4.1(iv)]{ZhanLoop}. By the same theorem, for $\Lcal$-almost every loop $\ell$, the trace $\ell\cap\CC$ has a $d$-dimensional Minkowski content measure, which we denote by $\mu_\ell$. Let $\Q_\kappa$ be the law of the curve $\eta$ in Theorem~\ref{thm:sssi}, and let $\Q_{\kappa,z}$ be its translate by $z$. On curves modulo reparametrization, the Campbell identity \cite[(4.4)]{ZhanLoop} reads
\begin{equation}\label{eq:Campbell}
 \Lcal(d\ell)\,\mu_\ell(dz)
 =\mathbb Q_{\kappa,z}(d\eta)\,dA(z).
\end{equation}

\begin{lemma}\label{lem:loop-finite}
There is a constant $C_*\in (0,\infty)$ depending only on $\kappa$ such that, for every $R>0$,
\[
 \Lcal\{\ell:\ell\cap\overline{B(0,R)}\ne\varnothing\}\le C_* R^{2-d}.
\]
\end{lemma}

\begin{proof}
Let $\mu_0^1$ denote the $\mathrm{SLE}_\kappa$ loop measure rooted at $0$ from \cite[Theorem~4.1(i)]{ZhanLoop}, and let $J(z)=-1/z$. By \cite[Theorem~4.1(iv)]{ZhanLoop}, $\Lcal=J(\mu_0^1)$. If a loop rooted at infinity hits $\overline{B(0,R)}$, its image under $J$ is rooted at $0$ and reaches $\partial B(0,1/R)$, so its diameter exceeds $1/(2R)$. The bound follows from \cite[Theorem~4.1(vi)]{ZhanLoop}.
\end{proof}

\begin{proposition} \label{prop:marked-ball}
Let $\eta$ be the sssi $\mathrm{SLE}_\kappa$ curve of Theorem~\ref{thm:sssi}, and let $\mu_\eta$ denote its occupation measure. There are constants $c_0,C_0\in(0,\infty)$ depending only on $\kappa$ such that, for every $r>0$,
\begin{equation}\label{eq:marked-ball}
 \E_{\Q_\kappa}\exp\left\{c_0\left(\frac{\mu_\eta(B(0,r))}{r^d}\right)^q\right\}\le C_0.
\end{equation}
\end{proposition}

\begin{proof}
For $R>0$, let $\Lcal_R$ be the restriction of $\Lcal$ to loops that hit $\partial B(0,R)$. For each such loop $\ell$, let $K_R(\ell)$ be its oriented initial segment from $\infty$ up to and including its first hit of $\partial B(0,R)$, and set $\Pi_R=(K_R)_*\Lcal_R$. By Lemma~\ref{lem:loop-finite}, $|\Pi_R|=|\Lcal_R|\le C_*R^{2-d}$. Applying inversion to the conformal Markov property in \cite[Theorem~4.1(ii)]{ZhanLoop} gives
\begin{equation}\label{eq:entrance-kernel}
 \Lcal_R
 =\int \Pi_R(d\alpha)\,
 \bigl(\alpha\oplus\mathbb P^{\#}_{U_\alpha;\alpha_{\rm tip}\to\infty}\bigr).
\end{equation}
Here $U_\alpha$ is the component of $\widehat{\CC}\setminus\alpha$ in which the tip of $\alpha$ and $\infty$ determine the marked prime ends, $\mathbb P^{\#}_{U_\alpha;\alpha_{\rm tip}\to\infty}$ is the law of chordal $\mathrm{SLE}_\kappa$ between these prime ends, and $\oplus$ denotes concatenation.

Take $R=8r$ and set $Y_\ell=\mu_\ell(B(0,2r))/(2r)^d$. The initial segment $\alpha$ is disjoint from the open ball $B(0,R)$. Since this ball is connected, it either lies in $U_\alpha$ or is disjoint from $U_\alpha$. In the latter case the chordal future cannot enter $B(0,R)$, so $Y_\ell=0$. In the former case the entire loop and its chordal future have the same trace in $B(0,R)$, so their Minkowski content measures agree there by restriction and uniqueness. Since $R=4(2r)$, Lemma~\ref{lem:interior-transport} bounds the expectation of $e^{c_1Y_\ell^q}$ over the chordal future by $C_1$, where $c_1>0$ and $C_1\ge1$ depend only on $\kappa$. Thus \eqref{eq:entrance-kernel} gives
\[
 \int e^{c_1Y_\ell^q}\,\Lcal_R(d\ell)\le C_1|\Pi_R|\le C_1C_*R^{2-d}.
\]

Take $\lambda=c_1/2$. Every loop visiting $B(0,r)$ hits $\partial B(0,R)$. Since $\mu_\ell(B(0,r))\le (2r)^dY_\ell$ and $y e^{\lambda y^q}\lesssim_\kappa e^{c_1y^q}$ for $y\ge0$, we obtain
\begin{equation}\label{eq:loop-weighted}
 \int \mu_\ell(B(0,r))e^{\lambda Y_\ell^q}\,\Lcal(d\ell)
 \lesssim_\kappa r^d\int e^{c_1Y_\ell^q}\,\Lcal_R(d\ell)
 \lesssim_\kappa r^d R^{2-d}\lesssim_\kappa r^2.
\end{equation}

Apply \eqref{eq:Campbell} to $\ind_{\{z\in B(0,r)\}}e^{\lambda Y_\ell^q}$. By \eqref{eq:loop-weighted}, the left-hand side is at most $C_\kappa r^2$. Since $\Q_{\kappa,z}$ is the translate of $\Q_\kappa$ by $z$ and $B(0,r)\subset B(0,2r)-z$ for every $z\in B(0,r)$, the right-hand side satisfies
\begin{equation}\label{eq:campbell-marked}
 \int_{B(0,r)} \E_{\Q_\kappa}\exp\left\{\lambda\left(\frac{\mu_\eta(B(0,2r)-z)}{(2r)^d}\right)^q\right\}\,dA(z)
 \ge \pi r^2\,\E_{\Q_\kappa}\exp\left\{2^{-dq}\lambda\left(\frac{\mu_\eta(B(0,r))}{r^d}\right)^q\right\}.
\end{equation}
Dividing \eqref{eq:campbell-marked} by $\pi r^2$ proves \eqref{eq:marked-ball} with $c_0=2^{-dq}\lambda$.
\end{proof}

\subsection{Stationary mass transport}\label{sec:sssi-mass-transport}
For a continuous path $x:\R\to\CC$ and $t\in\R$, define its recentered time shift by
\[
 (\theta_t x)(s)=x(t+s)-x(t),\qquad s\in\R.
\]
If $x(0)=0$, then
\begin{equation}\label{eq:backpoint}
 -(\theta_t x)(-t)=x(t).
\end{equation}

\begin{lemma}\label{lem:mass-transport}
Let $X:\R\to\CC$ be a continuous random curve with $X(0)=0$ and $\theta_tX\stackrel{\rm law}=X$ for every $t\in\R$. For $N>0$, let $\nu_N$ be the pushforward under $X$ of Lebesgue measure on $[-N,N]$. Then for every nonnegative measurable functional $F$ on $C(\R,\CC)$ and every deterministic Borel set $D\subset\CC$,
\begin{equation}\label{eq:mass-transport}
 \E\int_{-N}^N \ind_{\{X(t)\in D\}}F(\theta_tX)\,dt
 =\E\bigl[F(X)\nu_N(-D)\bigr],
\end{equation}
where $-D=\{-z:z\in D\}$.
\end{lemma}

\begin{proof}
For each $t\in\R$, \eqref{eq:backpoint} and $\theta_tX\stackrel{\rm law}=X$ give
\[
 \E\left[\ind_{\{X(t)\in D\}}F(\theta_tX)\right]
 =\E\left[\ind_{\{-(\theta_tX)(-t)\in D\}}F(\theta_tX)\right]
 =\E\left[\ind_{\{X(-t)\in -D\}}F(X)\right].
\]
Integrating over $[-N,N]$, applying Tonelli's theorem, and changing variables from $t$ to $-t$ give \eqref{eq:mass-transport}.
\end{proof}

\begin{proposition}\label{prop:ball-to-set}
Let $X:\R\to\CC$ be a continuous random curve with $X(0)=0$ and
$\theta_tX\stackrel{\rm law}=X$ for every $t\in\R$, and let $\nu$ be its
occupation measure.  Suppose that, for some $d,q,c_0,C_0>0$,
\begin{equation}\label{eq:abstract-ball}
 \E\exp\left\{c_0
 \left(\frac{\nu(B(0,r))}{r^d}\right)^q\right\}\le C_0,
 \qquad r>0.
\end{equation}
Then there are $c,C>0$, depending only on $c_0,C_0,d,q$, such that for every
deterministic bounded Borel set $D\subset\CC$ of diameter $L>0$,
\begin{equation}\label{eq:abstract-set}
 \E\exp\left\{c
 \left(\frac{\nu(D)}{L^d}\right)^q\right\}\le C.
\end{equation}
\end{proposition}

\begin{proof}
Let $\nu_N$ denote the occupation measure of $X$ on $[-N,N]$, and set
\[
 M_N=\nu_N(D),\qquad \widetilde M_N=\nu_N(-D),\qquad Z_L=\nu(B(0,2L)).
\]
If $X(t)\in D$, then $D\subset B(X(t),2L)$, so
\begin{equation}\label{eq:pathwise-domination}
 M_N\le \nu(B(X(t),2L)).
\end{equation}
Fix $s>1$. For a path $x$, let $F(x)$ be its occupation time in $B(0,2L)$ raised to the power $s-1$. Then $F(X)=Z_L^{s-1}$ and $F(\theta_tX)=\nu(B(X(t),2L))^{s-1}$. By \eqref{eq:pathwise-domination}, Lemma~\ref{lem:mass-transport}, and H\"older's inequality,
\begin{align*}
 \E[M_N^s]
 &=\E\int_{-N}^N \ind_{\{X(t)\in D\}}M_N^{s-1}\,dt\le \E\int_{-N}^N \ind_{\{X(t)\in D\}}F(\theta_tX)\,dt\\
 &
 =\E[Z_L^{s-1}\widetilde M_N]
 \le \|Z_L\|_s^{s-1}\|\widetilde M_N\|_s.
\end{align*}
Applying the same estimate to $-D$ interchanges $M_N$ and $\widetilde M_N$. Thus $A_N:=\max\{\|M_N\|_s,\|\widetilde M_N\|_s\}$ satisfies $A_N^s\le\|Z_L\|_s^{s-1}A_N$. Since $A_N\le 2N$, by \eqref{eq:abstract-ball}, we obtain $A_N\le\|Z_L\|_s$. Letting $N\uparrow\infty$ and using monotone convergence gives $\E[\nu(D)^s]\le \E[Z_L^s]$, $s>1$.

Set $U=c_0(\nu(D)/(2L)^d)^q$ and $V=c_0(Z_L/(2L)^d)^q$. By \eqref{eq:abstract-ball}, $\E[e^V]\le C_0$. For every positive integer $n$ with $nq>1$, the moment comparison with $s=nq$ gives $\E[U^n]\le\E[V^n]$. If $q>1$, summing the exponential series gives $\E[e^U]\le\E[e^V]\le C_0$. For $0<q\le1$, let $m$ be the least integer with $mq>1$. For $1\le n<m$, we have $U^n\le 1+U^m$, so $\E[U^n]\le 1+\E[V^m]\le 1+m!C_0$. Summing the series again gives $\E[e^U]\le C$ with $C$ depending only on $q$ and $C_0$. This proves \eqref{eq:abstract-set} with $c=2^{-dq}c_0$.
\end{proof}

\begin{proof}[Proof of Theorem~\ref{thm:sssi}]
By \eqref{eq:sssi-intro}, $\theta_t\eta$ has the same law as $\eta$ for every $t\in\R$. Proposition~\ref{prop:marked-ball} gives \eqref{eq:abstract-ball} for $X=\eta$ and $\nu=\mu_\eta$, so Proposition~\ref{prop:ball-to-set} yields \eqref{eq:sssi-main}.
\end{proof}

\section{Sharpness}\label{sec:sharpness}
Throughout this section, $\gamma$ denotes chordal $\mathrm{SLE}_\kappa$ in $\HH$ from $0$ to $\infty$, and $\mu_\gamma$ denotes the $d$-dimensional Minkowski content measure on $\gamma$.

\begin{proposition} \label{input:dense}
Fix $s_0\in(0,1]$. There are constants $a,b,C>0$ and an integer $m_0$, depending only on $\kappa$ and $s_0$, such that for every integer $m\ge m_0$ there is a deterministic $v_m>0$ with the following property. For every $w\in\HH$ with $|w|=1$ and $\Im w\ge s_0$, put $K_w=\overline{B(w,s_0/2)}$. Then, for every $0<v\le v_m$,
\BGE \label{eq:finite-mesh-dense}
 \PP\left\{v^{-p}\area(\gamma^v\cap K_w)\ge a m^p\right\}
 \ge C e^{-bm^2}.
\EDE
\end{proposition}

\begin{proof}
Fix $w\in\HH$ with $|w|=1$ and $y:=\Im w\ge s_0$. We use the dense-visit construction from the proof of \cite[Proposition~6.4]{OAI}. The proposition itself estimates the probability under the finite centered law of accumulating sausage area in a larger ball. The estimate we need is implicit in the proof: under the ordinary chordal law, the event constructed there has probability at least $c e^{-Cm^2}$, and the sausage area counted in \cite[(6.16)]{OAI} lies inside $K_w$. For the reader's convenience, we recall the construction and identify the strips inside $K_w$ that give the area lower bound.

By reflection symmetry of chordal $\mathrm{SLE}_\kappa$ about the imaginary axis, we may assume that $\Re w\ge0$. Set $x_0=\Re w+s_0/4$, and let $D$ be the closed disk of radius $s_0/20$ centered at $x_0+iy$. Since $s_0/4+s_0/20=3s_0/10<s_0/2$, we have $D\subset K_w$. Let $m$ be a positive integer and put $h_0=s_0/200$. For $j=0,\ldots,2m$, let $I_j$ be the vertical segment from $X_j+i(y-h_0)$ to $X_j+i(y+h_0)$, where $X_j=x_0+h_0(j/m-1)$.

To connect $0$ to the lower endpoint of $I_0$, set $R=h_0/2$ and $H=s_0/3$. Go vertically from $0$ to $i(H-R)$, then follow a quarter-circle of radius $R$ to $R+iH$, turning to the right. Continue horizontally to $X_0-R+iH$, then follow a second quarter-circle of radius $R$ to $X_0+i(H+R)$, turning upwards. Continue vertically to the lower endpoint of $I_0$ and traverse $I_0$ upwards. For $0\le j<2m$, join $I_j$ to $I_{j+1}$ by a semicircle of radius $h_0/(2m)$, above their upper endpoints when $j$ is even and below their lower endpoints when $j$ is odd. After traversing $I_{2m}$, continue vertically upwards to infinity.

Fix a sufficiently large absolute constant $C_0\ge20$, set $\ell=h_0/(C_0m)$, and rescale the plane by $z\mapsto z/\ell$. Write $\widehat\gamma=\ell^{-1}\gamma$, $\widehat K=\ell^{-1}K_w$, and $\widehat v=v/\ell$. By scale invariance, $\widehat\gamma$ has the same law as $\gamma$ as an unparametrized curve. Since $2-p=d$, we have
\[
 v^{-p}\area(\gamma^v\cap K_w)
 =\ell^d\,\widehat v^{-p}\area(\widehat\gamma^{\widehat v}\cap\widehat K).
\]
Since $\ell^d m^2=(h_0/C_0)^d m^p$, it suffices for \eqref{eq:finite-mesh-dense} to prove that there exist constants $c,C>0$, depending only on $\kappa$ and $s_0$, such that
\BGE \label{rescaled}
 \PP\left\{\widehat v^{-p}\area(\widehat\gamma^{\widehat v}\cap\widehat K)\ge c m^2\right\}\ge c e^{-Cm^2}
\EDE
for all sufficiently large $m$ and all sufficiently small $\widehat v>0$, with the mesh cutoff uniform in $w$ for each $m$.

Under this dilation, the ribbon of half-width $\ell$ about the deterministic curve becomes a ribbon of half-width one. Parametrize its centerline by arclength from $0$, writing $\sigma(s)$ for the point at arclength $s$ and $n(s)=i\sigma'(s)$ for the unit normal there. The normal coordinates on the rescaled ribbon are given by
\[
 P(s,u)=\sigma(s)+u n(s),\qquad s\ge0,\quad -1\le u\le1.
\]
After rescaling, adjacent vertical segments in the array are $C_0$ apart. The semicircles joining them have radius $C_0/2$, and the two quarter-circles in the connector have radius $R/\ell=C_0m/2$. The connector approaches the first segment from below, and the final ray extends upwards from the last segment. The ribbon consists of rectangles along the straight segments and annular sectors along the circular arcs. By the spacing of the segments and the choice of turn radii, only consecutive pieces meet, and their intersection is their common normal cross-section. Thus $P$ is injective. The centerline begins with a vertical segment of length $(H-R)/\ell>4$ and has curvature at most $2/C_0\le1/10$ almost everywhere. On each arclength interval of length four, its tangent turns by at most $8/C_0$, so the restrictions of $P$ to $[s,s+4]\times[-1,1]$ have a common bi-Lipschitz bound depending only on $C_0$. Thus the rescaled ribbon satisfies the geometric assumptions stated in the opening paragraph of \cite[Section~6.1, p.~32]{OAI}, with bounds uniform in $m$ and $w$.

In the rescaled coordinates, take the cross-sections $P(\{k\}\times[-1,1])$ at integer arclengths $k=0,\ldots,N$, where $N=\lceil s_{\rm out}\rceil+1$ and $s_{\rm out}$ is the arclength at which the final vertical ray exits $\ell^{-1}D$. In the original coordinates the connector has length $O(1)$, the $2m+1$ vertical segments have total length $2(2m+1)h_0$, and the semicircles and the final extension up to the disk exit have total length $O(h_0)$. Consequently $N=O((1+mh_0)/\ell)=O(m^2)$ uniformly in $w$. Each rescaled vertical segment has length $2C_0m$ and hence contains at least $2C_0m-2\ge C_0 m$ full unit steps. Since the original segments lie within distance $\sqrt{2}h_0$ of the center of $D$ and $D$ has radius $10h_0$, the ribbon along each rescaled segment lies in $\ell^{-1}D$. Thus at least $2C_0m^2$ steps have their middle strips contained in $\ell^{-1}D$.

Set $\epsilon=1/50$ and $\tau_0=0$. For $1\le k\le N$, let $\tau_k$ be the first time $\widehat\gamma$ hits the cross-section $s=k$. For $0\le k<N$, the step from $s=k$ to $s=k+1$ satisfies the passage conditions if $\tau_k<\tau_{k+1}<\infty$, the trace stays in the part of the ribbon with $k-1/4<s<k+1$ and $|u|<\epsilon$ throughout $[\tau_k,\tau_{k+1})$, and $\widehat\gamma(\tau_{k+1})=P(k+1,u)$ for some $|u|<\epsilon/4$. For $k=0$, these conditions alone define success. For $1\le k<N$, a successful step must also satisfy \BGE \widehat v^{-p}\area\bigl(\widehat\gamma_k^{\,\widehat v}\cap D_k\bigr)\ge c_1\label{step-lb}\EDE
where $\widehat\gamma_k:=\widehat\gamma[\tau_k,\tau_{k+1}]$ and $D_k:=P\bigl((k+1/3,k+2/3)\times(-\epsilon,\epsilon)\bigr)$.

By vertical-slit support, the initial step succeeds with probability at least $q_0>0$, uniformly in $m$ and $w$, since the rescaled ribbon is straight vertical up to $s=4$ \cite[proof of Proposition~6.4, p.~42]{OAI}. By \cite[Lemma~6.1]{OAI}, we may choose $c_1>0$ and decrease $q_0$ if necessary so that, for each $1\le k<N$, on the event that the preceding steps have succeeded, the conditional probability of success given the trace up to $\tau_k$ is at least $q_0$. This holds for all sufficiently small $\widehat v>0$, with the constants and the mesh cutoff uniform in $k,m,w$. By iteration, all $N$ steps succeed with probability at least $q_0^N$.

On this event, sum the area bounds over the indices $1\le k<N$ for which $D_k\subset\ell^{-1}D$, as in the derivation of \cite[(6.16)]{OAI}. Every step counted on the vertical segments of the array occurs after the initial step and hence has index $k\ge1$. Thus there are at least $2C_0m^2$ such indices. Since the strips $D_k$ are disjoint, $\ell^{-1}D\subset\widehat K$, and each $\widehat\gamma_k$ is a segment of $\widehat\gamma$, we obtain by \eqref{step-lb}
\[
 \widehat v^{-p}\area\bigl(\widehat\gamma^{\,\widehat v}\cap\widehat K\bigr)
 \ge \sum_{k:\,D_k\subset\ell^{-1}D}
       \widehat v^{-p}\area\bigl(\widehat\gamma_k^{\,\widehat v}\cap D_k\bigr)
 \ge 2C_0c_1m^2.
\]
Since $N\le C_*m^2$ uniformly in $w$, all $N$ steps succeed with probability at least $q_0^N\ge e^{-Cm^2}$ for some $C>0$. This proves the rescaled estimate \eqref{rescaled} and hence the original estimate \eqref{eq:finite-mesh-dense}.
\end{proof}

\begin{lemma} \label{lem:chordal-dense}
Fix $s_0\in(0,1]$. There exist constants $a,b,C>0$ and an integer $m_0\ge1$, depending only on $\kappa$ and $s_0$, such that for every $w\in\HH$ with $\sin(\arg w)\ge s_0$ and every $m\ge m_0$,
\[
 \PP\{\mu_\gamma(K_w)\ge a|w|^d m^p\}\ge C e^{-b m^2},
 \qquad K_w:=\overline{B(w,s_0|w|/2)}.
\]
\end{lemma}

\begin{proof}
Put $r=|w|$ and $\widetilde\gamma=r^{-1}\gamma$. Since $r^{-1}K_w=K_{w/r}$, the conformal covariance \eqref{eq:interior-mc} gives $\mu_{\widetilde\gamma}(K_{w/r})=r^{-d}\mu_\gamma(K_w)$. By the scaling property of chordal $\mathrm{SLE}_\kappa$, $\widetilde\gamma$ has the same law as $\gamma$ up to reparametrization. Thus $r^{-d}\mu_\gamma(K_w)$ has the same law as $\mu_\gamma(K_{w/r})$, and it suffices to consider $|w|=1$. Fix such a $w$ with $\Im w\ge s_0$, and put $K'_w=\overline{B(w,s_0/4)}$. Proposition~\ref{input:dense}, applied with $s_0/2$ in place of $s_0$, gives \eqref{eq:finite-mesh-dense} with $K'_w$ in place of $K_w$.

For $0<v<s_0/4$, we have $\gamma^v\cap K'_w\subset(\gamma\cap K_w)^v$. By \eqref{lin-Cont},
\begin{equation}\label{bound}
 \limsup_{v\downarrow0}v^{-p}\area(\gamma^v\cap K'_w)
 \le \overline{\Cont}_d(\gamma\cap K_w)
 \le \mu_\gamma(K_w).
\end{equation}
Fix an integer $m\ge m_0$ and let $v_m$ be as in Proposition~\ref{input:dense}. For $n\ge1$, set $\varepsilon_n=2^{-n}v_m$ and let $E_n=\{\varepsilon_n^{-p}\area(\gamma^{\varepsilon_n}\cap K'_w)\ge a m^p\}$. By Proposition~\ref{input:dense} with parameter $s_0/2$, $\PP(E_n)\ge C e^{-bm^2}$ for every $n$. Since the unions $\bigcup_{n\ge j}E_n$ decrease to $\limsup_{n\to\infty}E_n$ and each contains $E_j$, continuity from above gives $ \PP\bigl(\limsup E_n\bigr)
 =\lim_{j\to\infty}\PP\Bigl(\bigcup_{n\ge j}E_n\Bigr)
 \ge C e^{-bm^2}$.

On $\limsup E_n$, the normalized area is at least $a m^p$ for infinitely many $n$, so \eqref{bound} gives $\mu_\gamma(K_w)\ge a m^p$ almost surely on this event. This proves the assertion for integer $m$. For real $m\ge m_0$, apply the integer bound at $\lceil m\rceil$ and use $m\le\lceil m\rceil\le2m$, replacing $b$ by $4b$.
\end{proof}

We next extend Lemma~\ref{lem:chordal-dense} to disks compactly contained in simply connected domains, with explicit dependence on the ratio of the disk radius to the conformal radius at its center and on the angle factor.

\begin{proposition}\label{prop:chordal-ball-lower}
Let $U\subsetneq\CC$ be simply connected with distinct boundary prime ends $u,v$. Let $\gamma_U$ be chordal $\mathrm{SLE}_\kappa$ in $U$ from $u$ to $v$, and let $\mu_U$ be the $d$-dimensional Minkowski content measure on $\gamma_U\cap U$. There exist $a,b,C>0$ and an integer $m_0\ge1$, depending only on $\kappa$, such that for every $z\in U$, every $0<r<\dist(z,\partial U)$, and every integer $m\ge m_0$,
\begin{equation}\label{eq:chordal-domain-lower}
 \PP\left\{r^{-d}\mu_U(B(z,r))\ge a m^p\right\}
 \ge C\left(\frac{r}{\crad_U(z)}\right)^p
 S_U(z;u,v)^{8/\kappa-1}e^{-bm^2},
\end{equation}
where $\crad_U$ and $S_U$ are defined just after \eqref{eq:green}.
\end{proposition}

\begin{proof}
Write $h=\crad_U(z)$, $\rho=r/2$, and $\beta=8/\kappa-1$. Let $U_s$ be the remaining domain of $\gamma_U$ at time $s$, and let $\tau$ be the first time at which $z\in U_s$ and $\crad_{U_s}(z)=\rho$, with $\tau=\infty$ if there is no such time. Put $E=\{\tau<\infty,\ S_{U_\tau}(z;\gamma_U(\tau),v)\ge1/2\}$.

Let $\PP_z^*$ denote the law of two-sided radial $\mathrm{SLE}_\kappa$ in $U$ from $u$ to $v$ through $z$. Under $\PP_z^*$, $\tau<\infty$ almost surely. After mapping $U$ conformally onto the unit disk with $z$ sent to $0$, the conformal radius of the remaining domain at $0$ at time $\tau$ is $\rho/h<1/2$, since $\rho=r/2$ and $h\ge\dist(z,\partial U)>r$. The density comparison immediately following \cite[(40)]{LR}, integrated over angles in $[\pi/6,5\pi/6]$, therefore gives $\PP_z^*(E)\ge c_\kappa$ uniformly in the initial angle.

By the Green's function change of measure in \cite[(38)]{LR} and the formula \eqref{eq:green},
\[
 \PP_z^*(E)
 =\int_E \frac{G_{U_\tau}(z;\gamma_U(\tau),v)}{G_U(z;u,v)}\,d\PP
 \le \left(\frac{h}{\rho}\right)^p S_U(z;u,v)^{-\beta}\PP(E),
\]
where we used $S_{U_\tau}(z;\gamma_U(\tau),v)\le1$. Combining the two bounds and substituting $\rho=r/2$ gives
\BGE\label{eq:E-lower}
 \PP(E)\ge 2^{-p}c_\kappa\left(\frac{r}{h}\right)^p S_U(z;u,v)^\beta.
\EDE

On $E$, conditional on $\FF_\tau$, the future is chordal $\mathrm{SLE}_\kappa$ in $U_\tau$ from $\gamma_U(\tau)$ to $v$. Choose a conformal map $f:\HH\to U_\tau$ taking $0,\infty$ to these prime ends, respectively, and normalize it so that $w=f^{-1}(z)$ satisfies $|w|=1$. Put $y=\Im w$. Then $1/2\le y\le1$ and $2y|f'(w)|=\rho$. The disk $K_w=\overline{B(w,1/4)}$ in Lemma~\ref{lem:chordal-dense} with $s_0=1/2$ is contained in $\overline{B(w,y/2)}$. Applying the Koebe growth and distortion estimates to $f$ on $B(w,y)$ gives, for every $\zeta\in K_w$,
\[
 |f(\zeta)-z|\le2y|f'(w)|=\rho<r,
 \qquad
 |f'(\zeta)|\ge\frac{4}{27}|f'(w)|\ge\frac{r}{27}.
\]
Since $f(K_w)\subset B(z,r)$, conformal covariance \eqref{eq:interior-mc} gives
\[
 \mu_U(B(z,r))\ge\int_{K_w}|f'(\zeta)|^d\,\mu_\gamma(d\zeta)
 \ge (r/27)^d\mu_\gamma(K_w).
\]
Lemma~\ref{lem:chordal-dense} therefore gives constants $a,b,C_1\in (0,\infty)$ depending only on $\kappa$, such that on $E$,
\[
 \PP\!\left\{r^{-d}\mu_U(B(z,r))\ge a m^p\,\middle|\,\FF_\tau\right\}
 \ge C_1e^{-bm^2},\qquad m\ge m_0,
\]
Integrating over $E$ and using \eqref{eq:E-lower} proves \eqref{eq:chordal-domain-lower}.
\end{proof}

\begin{proposition}\label{prop:sssi-ball-lower}
For sssi $\mathrm{SLE}_\kappa$, there exist $a,b,C>0$ and an integer $m_0\ge1$, depending only on $\kappa$, such that \eqref{eq:sharp-sssi} holds for every $z\in\CC$, every $r>0$, and every integer $m\ge m_0$.
\end{proposition}

\begin{proof}
We first consider $0<r\le|z|$. By \eqref{eq:sssi-intro}, $r^{-d}\mu_\eta(B(z,r))$ has the same law as $t^{-d}\mu_\eta(B(1,t))$, where $t=r/|z|\le1$. We will show that
\[
 \PP\{t^{-d}\mu_\eta(B(1,t))\ge a m^p\}\ge C t^p e^{-bm^2},
 \qquad 0<t\le1,\quad m\ge m_0.
\]

Let $\alpha$ be the first arm of $\eta$, from $\infty$ to $0$. Given $\alpha$, the second arm is chordal $\mathrm{SLE}_\kappa$ from $0$ to $\infty$ in the remaining domain $U$ \cite[Section~2.2]{ZhanLoop}. We claim that $\PP\{1\in U\}>0$. Indeed, rotation and scaling invariance imply that $\PP\{w\in U\}$ is the same for every $w\ne0$. If this probability were zero, then almost surely $U$ would contain no point of a fixed countable dense subset of $\CC\setminus\{0\}$, whereas $U$ is a nonempty open set.

For $j\ge2$, let
\[
 G_j=\left\{1\in U,\ \dist(1,\partial U)>j^{-1},\
 S_U(1;0,\infty)>j^{-1}\right\}.
\]
These events increase to $\{1\in U\}$, so we may fix $j$ with $\PP(G_j)>0$. Since $0\in\partial U$, the Koebe quarter theorem gives $\crad_U(1)\le4\dist(1,\partial U)\le4$ on $G_j$. For every $0<t\le1$, we have $t/j<\dist(1,\partial U)$ on $G_j$, so Proposition~\ref{prop:chordal-ball-lower} applies to the second arm in $B(1,t/j)$. Since $B(1,t/j)\subset B(1,t)$, we obtain, with $\beta=8/\kappa-1$,
\[
 \PP\!\left\{\mu_\eta(B(1,t))\ge a_1(t/j)^d m^p
       \,\middle|\,\alpha\right\}
 \ge C_1\left(\frac{t}{4j}\right)^p j^{-\beta}e^{-b_1m^2}
 \qquad\text{on }G_j.
\]
Integrating over $G_j$ gives the desired estimate at $1$, since $j$ and $\PP(G_j)$ depend only on $\kappa$. Scaling back with $t=r/|z|$ gives \eqref{eq:sharp-sssi} for $0<r\le|z|$.

It remains to consider $|z|<r$. Then $1\wedge \frac r{|z|}=1$. Choose $u$ on the ray from $0$ through $z$ with $|u|=r/2$, taking $u=r/2$ when $z=0$. Then $|u-z|\le r/2$, so $B(u,r/4)\subset B(z,r)$. Applying the preceding case to $B(u,r/4)$ gives
\[
 \PP\{\mu_\eta(B(z,r))\ge a(r/4)^d m^p\}
 \ge C\,2^{-p}e^{-bm^2}.
\]
Thus \eqref{eq:sharp-sssi} holds in both cases with the common constants $4^{-d}a$ and $2^{-p}C$ in place of $a$ and $C$, respectively, and with the same $b$ and $m_0$.
\end{proof}

\begin{proof}[Proof of Theorem~\ref{thm:sharpness}]
We have $\crad_{\HH}(z)=2\Im z$ and $S_{\HH}(z;0,\infty)=\Im z/|z|$. Thus Proposition~\ref{prop:chordal-ball-lower} gives \eqref{eq:sharp-chordal}, with the factor $2^{-p}$ absorbed into $C$. Proposition~\ref{prop:sssi-ball-lower} gives \eqref{eq:sharp-sssi}. We may use common constants in the two estimates by taking the smaller values of $a,C$ and the larger values of $b,m_0$.

Fix a disk $B(z,r)$ in either statement and put $X=r^{-d}\mu(B(z,r))$. For this disk, the corresponding lower bound has the form $\PP\{X\ge am^p\}\ge c(z,r)e^{-bm^2}$ with $c(z,r)>0$. By \eqref{eq:dq-intro}, $pq=2$, so for every integer $m\ge m_0$,
\BGE \label{contra}
 \E\!\left[\exp\{b(X/a)^q\}\ind_{\{X\ge am^p\}}\right]
 \ge e^{bm^2}\PP\{X\ge am^p\}\ge c(z,r).
\EDE
 If $\E\exp\{b(X/a)^q\}$ were finite, then $X$ would be finite almost surely and the left-hand side of \eqref{contra} would tend to zero as $m\to\infty$ by dominated convergence, contradicting the lower bound $c(z,r)>0$. Thus $\E\exp\{b(X/a)^q\}=\infty$ for every such disk. Given a nonempty deterministic open set $D$ in the corresponding domain, choose a disk $B(z,r)$ with $\overline{B(z,r)}\subset D$. Since $\mu(D)\ge\mu(B(z,r))$, we obtain \eqref{eq:critical-coefficient-finite} with $c_D=b/(a r^d)^q$.
\end{proof}

\section{Uniform local occupation bounds}\label{sec:uniform-local}

\begin{proposition} \label{prop:uniform-local-upper}
Let $\mu$ denote either $\mu_\gamma$ or $\mu_\eta$. There exists a deterministic constant $C_\kappa>0$ such that, almost surely, for every $R>0$, there is a random radius $r_R\in (0,1)$ such that
\[
 \sup_{z\in B(0,R)}\mu(B(z,r))
 \le C_\kappa r^d\bigl(\log(1/r)\bigr)^{1/q},\quad 0<r<r_R.
\]
\end{proposition}

\begin{proof}
Since $\diam B(w,3r)=6r$, Theorems~\ref{thm:main} and~\ref{thm:sssi} and Markov's inequality give constants $c,C>0$, depending only on $\kappa$, such that
\begin{equation}\label{eq:ball-tail-grid}
 \PP\{\mu(B(w,3r))>u r^d\}\le C e^{-cu^q},
 \qquad w\in\CC,\quad r>0,\quad u\ge0.
\end{equation}
In the chordal case we apply Theorem~\ref{thm:main} to $B(w,3r)\cap\HH$ whenever this intersection is nonempty.

Fix  $R\in\N$,  and let $G_n=2^{-n}(\Z+i\Z)\cap B(0,R+1)$, $n\in\N$. Then $|G_n|\lesssim_R2^{2n}$, and every $z\in B(0,R)$ lies within distance $2^{-n}$ of some $w\in G_n$. Choose $A>0$ so that $cA^q>2\log 2$. By \eqref{eq:ball-tail-grid} and the union bound,
\[
 \sum_{n=1}^\infty
 \PP\left\{\max_{w\in G_n}\mu(B(w,3\cdot 2^{-n}))>A 2^{-nd} n^{1/q}\right\}
 \lesssim_R\sum_{n=1}^\infty 2^{2n}e^{-cA^q n}<\infty.
\]
For each $z\in B(0,R)$, choose $w\in G_n$ with $|z-w|\le 2^{-n}$. Then $B(z,2^{-n})\subset B(w,3\cdot 2^{-n})$, so the Borel--Cantelli lemma gives, almost surely, a random integer $N_R\ge1$ such that
\[
 \sup_{z\in B(0,R)}\mu(B(z,2^{-n}))
 \le A 2^{-nd}n^{1/q},\qquad n\ge N_R.
\]
Set $r_R=2^{-N_R}$. For $0<r<r_R$, choose $n$ so that $2^{-n-1}<r\le 2^{-n}$. Then $n\ge N_R$, $2^{-n}<2r$, and $n\le\log(1/r)/\log 2$. By monotonicity,
\[
 \sup_{z\in B(0,R)}\mu(B(z,r))
 \le A 2^{-nd}n^{1/q}
 \le C_\kappa r^d\bigl(\log(1/r)\bigr)^{1/q},
\]
where $C_\kappa=2^d A(\log 2)^{-1/q}$. Taking the intersection of these probability-one events over all positive integers $R$, and using $B(0,R)\subseteq B(0,\lceil R\rceil)$, gives the assertion simultaneously for every $R>0$.
\end{proof}

\begin{corollary}[Two-sided local occupation bounds]\label{cor:natural-local-bounds}
Let $\xi$ be either the chordal curve $\gamma$ or the sssi curve $\eta$, parametrized by Minkowski content. For a nondegenerate compact interval $I$ in its natural-time domain, write $\mu_I^\xi $ for the pushforward of the Lebesgue measure on $I$ under $\xi$. In each of the following cases, the estimate holds almost surely for every such $I$.

\textup{(i)} For $\xi=\gamma$ and every $\varepsilon\in(0,d)$, there exist $c_\varepsilon>0$ and $r_0\in(0,1)$ such that
\begin{equation}\label{eq:time-local-upper}
 c_\varepsilon r^{d+\varepsilon}
 \le \mu_I^\gamma(B(z,r))
 \le C_\kappa r^d\bigl(\log(1/r)\bigr)^{1-d/2},
 \qquad z\in\gamma(I),\quad 0<r<r_0.
\end{equation}

\textup{(ii)} For $\xi=\eta$, there exist $c>0$ and $r_0\in(0,1)$ such that
\[
 c r^d\bigl(\log(1/r)\bigr)^{1-d}
 \le \mu_I^\eta(B(z,r))
 \le C_\kappa r^d\bigl(\log(1/r)\bigr)^{1-d/2},
 \qquad z\in\eta(I),\quad 0<r<r_0.
\]
The lower-bound constants and the radius $r_0$ may depend on
the sample path and on $I$, and in part \textup{(i)} also on
$\varepsilon$. The upper-bound constant $C_\kappa$ is deterministic
and depends only on $\kappa$.
\end{corollary}

\begin{proof}
Since $\xi(I)$ is bounded and $\mu_I^\xi\le\mu_\xi$, Proposition~\ref{prop:uniform-local-upper} gives both upper bounds, with $1/q=(2-d)/2$ by \eqref{eq:dq-intro}. The chordal lower bound follows from local H\"older continuity of every order below $1/d$, as shown in \cite[Proposition~4.8]{ZhanUniformDimension}. The sssi lower bound is given in \cite[Remark~7.2]{ZhanUniformDimension} as a consequence of the modulus of continuity in \cite[Theorem~1.2]{HoldenYuan}.
\end{proof}

\begin{remark}[Conformal transfer of the chordal bounds]\label{rem:local-upper-conformal}
Let $F:\HH\to D$ be conformal, and let $I$ be a nondegenerate compact natural-time interval such that $\gamma(I)\Subset\HH$. By \eqref{eq:interior-mc}, the occupation measure of the image segment parametrized by Minkowski content is $\nu=F_*(|F'|^d\mu_I^\gamma)$. Choose $L\ge1$ such that $F$ is $L$-bi-Lipschitz on a neighborhood of $\gamma(I)$. Then, for every $z\in\gamma(I)$,
\[
 L^{-d}\mu_I^\gamma(B(z,r/L))
 \le \nu(B(F(z),r))
 \le L^d\mu_I^\gamma(B(z,Lr)).
\]
Thus both bounds in \eqref{eq:time-local-upper} hold for the image segment, with constants depending on the segment and $F$.
\end{remark}

\subsection*{Use of artificial intelligence.}
The author used GPT~6 to explore proof ideas, identify potentially
relevant references, and assist with the exposition.
The author independently verified the references, simplified and
reorganized the arguments, and carefully checked all mathematical
details and the final manuscript.
The author takes full responsibility for its contents.

\end{document}